\documentclass[11pt,leqno]{article}
\usepackage{amsfonts}
\usepackage{amsmath}
\usepackage{amssymb}
\newcommand{\pref}[1]{(\ref{#1})}
\newtheorem{theo}{Theorem}[section]

\title{\Large \bf {The first eigenvalue of the Dirac operator on locally reducible Riemannian
manifolds}}
\author{{\sc Bogdan Alexandrov}
}
\date{}
\begin{document}
\maketitle
\vspace{5mm}
\begin{abstract}
We prove a lower estimate for the first eigenvalue of the Dirac operator on a compact locally reducible
Riemannian spin manifold with positive scalar curvature. We determine also the universal covers of the 
manifolds on which the smallest possible eigenvalue is
attained.
\\[10mm]
\\[10mm]
\end{abstract}

\section{Introduction}

Let $M$ be a compact Riemannian spin manifold with positive scalar curvature $s$. A well known
consequence of the Schr\"odinger--Lichnerowicz formula \cite{S,L}
\begin{equation}\label{0}
D^2 = \nabla^* \nabla + \frac{s}{4}
\end{equation}
is the Lichnerowicz vanishing theorem \cite{L}: the kernel of the Dirac operator $D$ coincides with 
the space of parallel spinors. If this
kernel is non-trivial, then the scalar curvature must be zero and even more, $M$ must be Ricci-flat. 
The simply connected irreducible manifolds which admit parallel spinors were described by Wang 
\cite{W}. They are exactly the manifolds whose holonomy group is one of $SU(n)$, $Sp(n)$, $G_2$, 
$Spin(7)$. For each of these cases Wang found the dimension of the space of parallel spinors.

If the scalar curvature is not identically zero, then the kernel of $D$ is trivial. This holds true, 
in particular, if $s>0$, which we shall assume from now on. In this case an estimate for the square of
the first eigenvalue $\lambda$ of the Dirac operator was obtained by Friedrich \cite{F}:
\begin{equation}\label{1}
\lambda ^2 \ge \frac{n}{n-1} \cdot \frac{\min s}{4},
\end{equation}
where $n$ is the dimension of the manifold. The limiting case of this estimate (i.e., the case of 
equality in \pref{1}) is characterized by the existence of a non-trivial real Killing spinor on $M$. 
It was shown by B\"ar \cite{B} that the last property is equivalent to the existence of a
parallel spinor on the cone over $M$. Thus, using Wang's results, he classified the simply connected
manifolds admitting real Killing spinors. They are the sphere $S^n$, the Einstein-Sasakian manifolds, 
the 3-Sasakian manifolds, the 6-dimensional nearly K\"ahler manifolds and the 7-dimensional manifolds 
with nearly parallel vector cross product \cite{B}.

The manifolds with non-trivial real Killing spinors have constant positive scalar curvature. They are
furthermore Einstein, locally irreducible and their holonomy is $SO(n)$. In particular, they do not 
admit any parallel $k$-form for $k \not = 0,n$ \cite{H}. This shows that the estimate \pref{1} cannot
be sharp, for example, on K\"ahler or quaternionic K\"ahler manifolds. Indeed, better estimates have 
been proved in these cases by Kirchberg \cite{K1,K2} and Kramer, Semmelmann, Weingart \cite{KSW} 
respectively. 

Another situation in which \pref{1} is not sharp was considered by the author, G. Gran\-tcha\-rov and
S. Ivanov. In \cite{AGI} they showed that if $M$ admits a non-trivial parallel 1-form, then
\begin{equation}\label{2}
\lambda ^2 \ge \frac{n-1}{n-2} \cdot \frac{\min s}{4}
\end{equation}
and the universal cover of a manifold, on which the first eigenvalue is the smallest possible, is a 
Riemannian product
$\mathbb{R} \times N$, where $N$ is a simply connected manifold with a real Killing spinor. (In fact,
it was proved by Moroianu and Ornea \cite{MO} that the same holds true under the weaker assumption
that the 1-form is harmonic and has constant length.)

A manifold with a non-trivial parallel 1-form is locally a product of a 1-dimensional and an 
$(n-1)$-dimensional manifold. In \cite{Kim} Kim considered the more general situation where the
manifold is locally a product of Riemannian manifolds of arbitrary dimensions $n_1$ and $n_2$. More 
precisely, he proved the following generalization of \pref{2}. Let the tangent bundle $TM = T_1
\oplus T_2$, where $T_1$ and $T_2$ are parallel, $\dim T_i = n_i$, $i=1,2$, and $n_1 \le n_2$. Then
the first eigenvalue of the Dirac operator satisfies
\begin{equation}\label{3}
\lambda ^2 \ge \frac{n_2}{n_2 - 1} \cdot \frac{\min s}{4}.
\end{equation}
Kim also gave examples of manifolds on which the equality in \pref{3} is attained. These are
global products $M = M_1 \times M_2$, where $\dim M_1 = n_1$, $\dim M_2 = n_2$, $M_2$ admits a real
Killing spinor and $M_1$ has a parallel spinor if $n_1 < n_2$ or either a parallel or a real Killing 
spinor if $n_1 = n_2$. In
\cite{Kim} Kim stated the question of whether these examples give the general form of the universal 
coverings of the manifolds on which the limiting case of \pref{3} occurs.

In the present paper we show that the answer of this question is positive. More generally, we consider 
manifolds which are locally a product of $k$ Riemannian manifolds and find an estimate for the first 
eigenvalue of the Dirac operator. We also describe the universal covers of the manifolds on which the 
limiting case of this estimate is attained. We summarize our results in the following theorem. 
\begin{theo}\label{th1}
Let $M$ be a compact Riemannian spin manifold with positive scalar curvature $s$. Let $TM = T_1 \oplus
\dots \oplus T_k$, where $T_i$ are parallel distributions of dimension  $n_i$, $i=1,\dots ,k$, and
$n_1 \le \dots \le n_k$. Then the first eigenvalue $\lambda$ of the Dirac operator satisfies
\begin{equation}\label{4}
\lambda ^2 \ge \frac{n_k}{n_k - 1} \cdot \frac{\min s}{4} .
\end{equation}
If the equality in \pref{4} is attained, then the universal cover $\widetilde{M}$ of $M$ is isometric 
to a
product $M_1 \times \dots \times M_k$, where $\dim M_i = n_i$, $M_k$ has a real Killing spinor and for
$i<k$ $M_i$ has a parallel spinor if $n_i<n_k$ or either a parallel spinor or a real Killing spinor if
$n_i=n_k$.
\end{theo}
Notice that the estimate \pref{4} contains \pref{1}, \pref{2} and \pref{3} as special cases. We remark 
also that $n_k > 1$ because of the positivity of the scalar curvature and that compact manifolds on
which \pref{4} is an equality do exist. For example, $M = M_1 \times \dots \times M_k$, where $M_i$,
$i=1,\dots ,p$, is a compact manifold admitting a parallel spinor (e.g.,
$M_i = S^1 \times \dots \times S^1$ ($n_i$ times)) and $M_j$,
$j=p+1,\dots ,k$, is a compact manifold of dimension $n_k$ admitting a real Killing spinor (e.g.,
$M_j = S^{n_k}$).

The proof of \pref{4} uses a typical technique for such situations. We introduce certain suitable 
twistor operator and prove a Weitzenb\"ock formula including
this operator and the Dirac operator, which in its turn implies \pref{4}. This simple proof is 
somewhat different than the proof of \pref{3} in \cite{Kim} but its advantage is that it allows us to 
handle easily the limiting case.

\section{Preliminaries}

Let $M$ be a compact Riemannian spin manifold such that the tangent bundle
\begin{equation}\label{4.5}
TM = T_1 \oplus \dots \oplus T_k,
\end{equation}
where $T_i$ are parallel and pairwise orthogonal distributions of dimension  $n_i$, $i=1,\dots ,k$. 

We denote the pointwise inner products by $<\cdot , \cdot >$ and the corresponding norms by $|\cdot|$.
The global inner products will be denoted by $(\cdot , \cdot )$,  
$(\cdot , \cdot ) = \int_M <\cdot , \cdot > dVol_M$, and the corresponding global norms by 
$||\cdot||$. 

By $e_{1,1},\dots , e_{1,n_1},\dots , e_{k,1},\dots , e_{k,n_k}$ we will denote an adapted local
orthonormal
frame, i.e., such that $e_{i,1},\dots , e_{i,n_i}$ spans $T_i$. The dual frame will be denoted by
$e^{1,1},\dots , e^{1,n_1},\dots ,$ $e^{k,1},\dots ,e^{k,n_k}$. Since $T_1,  \dots , T_k$ are parallel,
we can always assume that $\nabla e_{i,l} = 0$ at a fixed point $x$.

Another consequence of the parallelism of $T_i$ is that $R(X,Y)=0$ whenever $X \in T_i$, $Y \in T_j$
with $i \not = j$. 

Let $s_i$ be the "scalar curvature" of $T_i$, i.e.,
$$s_i = \sum_{l,m =1}^{n_i} <R(e_{i,l},e_{i,m})e_{i,m},e_{i,l}>.$$
Hence the scalar curvature $s$ of $M$ is $s= \sum_{i =1}^{k} s_i$.

We denote by $p_i$ the orthogonal projections from $TM$ onto $T_i$ and from $T^*M$ onto $T_i^*$.

Let $\Sigma M$ be the spinor bundle of $M$,
$\nabla : \Gamma (\Sigma M) \longrightarrow \Gamma (T^* M \otimes \Sigma M)$ the covariant derivative
of the Levi-Civita connection and $\mu: T^*M \otimes \Sigma M \longrightarrow \Sigma M$ the Clifford
multiplication (we write also $\alpha \cdot \psi$ instead of $\mu (\alpha \otimes \psi)$). Thus the
Dirac operator \\
$D: \Gamma (\Sigma M) \longrightarrow \Gamma (\Sigma M) \mbox{ is given by } D = \mu \circ
\nabla, \mbox{ i.e., } D\psi = \sum_{i =1}^{k} \sum_{l=1}^{n_i} e^{i,l} \cdot \nabla_{e_{i,l}} \psi.$

We define
$\nabla_i : \Gamma (\Sigma M) \longrightarrow \Gamma (T_i^* \otimes \Sigma M)$ by
${\nabla_i}_X \psi := \nabla_{p_i (X)} \psi$ and \\
$D_i: \Gamma (\Sigma M) \longrightarrow \Gamma (\Sigma M)$ by $D_i := \mu \circ \nabla_i$, i.e., 
$D_i \psi = \sum_{l=1}^{n_i} e^{i,l} \cdot {\nabla_i}_{e_{i,l}} \psi$. \\
So we have
$$\nabla = \sum_{i =1}^{k} \nabla_i, \qquad D= \sum_{i =1}^{k} D_i.$$

Next we list several formulae whose proofs are straightforward and quite similar to those of the 
corresponding results about the Dirac operator $D$.

Let $X \in \Gamma (TM)$ (resp. $\alpha \in \Gamma (T^*M)$) be orthogonal to $T_i$ (resp. $T_i^*$) and
$(\nabla X)_x =0$ (resp. $(\nabla \alpha)_x =0$) at a given point $x$. Then
\begin{equation}\label{5}
(\nabla_X (D_i \psi))_x = (D_i (\nabla_X \psi))_x, \qquad
(\alpha \cdot (D_i \psi))_x = -(D_i (\alpha \cdot \psi))_x.
\end{equation}
This implies, in particular,
\begin{equation}\label{6}
D_i D_j + D_j D_i  = 0 \quad \mbox{for } i \not = j
\end{equation}
and therefore
\begin{equation}\label{7}
D^2 = \sum_{i =1}^{k} D_i^2.
\end{equation}
We also have
\begin{equation}\label{8}
||D_i \psi ||^2 = (D_i^2 \psi , \psi)
\end{equation}
and the following "partial" Schr\"odinger--Lichnerowicz formula
\begin{equation}\label{9}
D_i^2 = \nabla_i^* \nabla_i + \frac{s_i}{4}.
\end{equation}

\section{Proof of Theorem~\ref{th1}}

Let $\pi : T^*M \otimes \Sigma M \longrightarrow T^*M \otimes \Sigma M$ be defined by
$$\pi(\alpha \otimes \psi ) = \alpha \otimes \psi + \sum_{i =1}^{k} \frac{1}{n_i}
\sum_{l =1}^{n_i } e^{i,l} \otimes (e^{i,l} \cdot p_i(\alpha) \cdot \psi ).$$
Clearly, this definition does not depend on the choice of the adapted frame $\{ e_{i,l} \}$ and
$\pi (T^*M \otimes \Sigma M) \subset \ker \mu$.

We introduce the following "adapted" twistor operator
$$P: \Gamma (\Sigma M) \longrightarrow \Gamma (T^*M \otimes \Sigma M): \qquad P:= \pi \circ \nabla .$$
In other words,
$$P\psi = \nabla \psi + \sum_{i =1}^{k} \frac{1}{n_i}
\sum_{l =1}^{n_i } e^{i,l} \otimes (e^{i,l} \cdot D_i \psi ).$$
We easily see
$$|P\psi|^2 = |\nabla \psi|^2 - \sum_{i =1}^{k} \frac{1}{n_i} |D_i \psi|^2  .$$
Using the Schr\"odinger--Lichnerowicz formula \pref{0} and \pref{8}, this implies
$$||P\psi||^2 = \left( \left( D^2 - \sum_{i =1}^{k} \frac{1}{n_i} D_i ^2 \right) \psi, \psi \right) - 
\left( \frac{s}{4}\psi, \psi \right).$$
By \pref{7} $D_k^2 = D^2 - \sum_{i =1}^{k-1} D_i ^2$. Hence
$$||P\psi||^2 = \left( 1 - \frac{1}{n_k} \right) (D^2\psi, \psi)
- \sum_{i =1}^{k-1} \left( \frac{1}{n_i} - \frac{1}{n_k} \right) (D_i ^2 \psi, \psi)
 - \left( \frac{s}{4}\psi, \psi \right)$$
and therefore
$$||D\psi||^2 = \frac{n_k}{n_k - 1}||P\psi||^2
+ \sum_{i =1}^{k-1} \frac{n_k - n_i}{n_k - 1}||D_i \psi||^2
+ \frac{n_k}{n_k - 1} \left( \frac{s}{4}\psi, \psi \right).$$
Since $n_k \ge n_i$ for each $i$, we obtain that the first eigenvalue $\lambda$ of $D$ satisfies
\pref{4}. This proves the first part of the theorem.

Suppose now that for the first eigenvalue $\lambda$ we have
$\lambda ^2 = \frac{n_k}{n_k - 1} \cdot \frac{\min s}{4} .$
Then for each eigenspinor $\psi$ of $D$ for $\lambda$ we have $P\psi = 0$, $D_i \psi = 0$ if
$n_i < n_k$ and $s$ is constant on the support of $\psi$. But $D^2$ is an elliptic operator
of second order and it follows from \cite{A} that each
eigenspinor of $D$ does not vanish on a dense open subset of $M$. Thus the scalar curvature $s$ is 
constant on the whole manifold and the parallelism
of $T_i$ implies that the same is true for $s_i$, $i=1,\dots,k$.

The equation $P\psi = 0$ is equivalent to
\begin{equation}\label{10}
\nabla_i \psi + \frac{1}{n_i} \sum_{l =1}^{n_i } e^{i,l} \otimes (e^{i,l} \cdot D_i \psi ) = 0, \qquad
i=1,\dots,k.
\end{equation}

If $D_i \psi = 0$ (in particular, if $n_i < n_k$), then
\begin{equation}\label{11}
\nabla_i \psi = 0.
\end{equation}

Let $D_i \psi \not = 0$. This means that $n_i = n_k$. We can write \pref{10} as
\begin{equation}\label{12}
{\nabla_i}_X \psi + \frac{1}{n_i} p_i (X^\flat) \cdot D_i \psi = 0, \qquad X \in TM,
\end{equation}
where $X^\flat$ is the 1-form corresponding to $X$ via the metric. Since $D\psi$ is also an
eigenspinor of $D$ for the eigenvalue $\lambda$, \pref{12} is satisfied by it too:
$${\nabla_i}_X (D\psi) + \frac{1}{n_i} p_i (X^\flat) \cdot D_i D\psi = 0, \qquad X \in TM.$$ 
Now, if $X \in \Gamma (TM)$ and $(\nabla X)_x = 0$, then \pref{5} and \pref{6} imply that at $x$
$${\nabla_i}_X (D_i\psi) + \frac{1}{n_i} p_i (X^\flat) \cdot D_i^2 \psi
+ \sum_{j \not = i} D_j ({\nabla_i}_X \psi + \frac{1}{n_i} p_i (X^\flat) \cdot D_i \psi)= 0.$$
This, together with \pref{12}, shows that
\begin{equation}\label{15}
{\nabla_i}_X (D_i \psi) + \frac{1}{n_i} p_i (X^\flat) \cdot D_i^2 \psi = 0, \qquad X \in TM.
\end{equation}
By a straightforward computation, similar to the case of the usual twistor spinors, \pref{12} and 
\pref{9} yield
\begin{equation}\label{16}
D_i^2 \psi =\frac{n_i}{n_i - 1} \cdot \frac{s_i}{4} \psi.
\end{equation}
In particular, by \pref{8}, $\lambda_i^2 :=\frac{n_i}{n_i - 1} \cdot \frac{s_i}{4} > 0$.

Let $\varphi_{i\pm} := \psi \pm \frac{1}{\lambda_i} D_i \psi$. By \pref{12}, \pref{15} and \pref{16}
we obtain
\begin{equation}\label{17}
{\nabla_i}_X \varphi_{i\pm} = \mp \frac{\lambda_i}{n_i} p_i (X^\flat) \cdot \varphi_{i\pm} , \qquad X 
\in TM.
\end{equation}
At least one of $\varphi_{i+}$ and $\varphi_{i-}$ is not identically zero because the same is true for
$\psi$. Let $x \in M$ be such that, for example, $\varphi_{i+} (x) \not = 0$. The universal cover
$\widetilde{M}$ of $M$ is a Riemannian product  $\widetilde{M} = M_1 \times \dots \times M_k$ because
of the decomposition \pref{4.5}. Let the point $\widetilde{x} = (x_1,\dots ,x_k)$ project on $x$.
Denote by $f_i: M_i \longrightarrow \widetilde{M}$ the inclusion
$f_i (y) = (x_1,\dots ,x_{i-1},y,x_{i+1}\dots ,x_k)$. Consider the pull-back 
$(q \circ f_i)^* \Sigma M$ of $\Sigma M$ to $M_i$, where $q:\widetilde{M} \longrightarrow M$ is the
projection. The bundle $(q \circ f_i)^* \Sigma M$ is a Clifford module on $M_i$ and is therefore a sum
of finitely many, say $r_i$, copies of the spinor bundle $\Sigma M_i$ of $M_i$. Also the Levi-Civita
connection on $\Sigma M$ pulls back to the Levi-Civita connection on 
$(q \circ f_i)^* \Sigma M = \oplus_{l=1}^{r_i} \Sigma M_i$. Thus \pref{17} implies that
$(q \circ f_i)^* \varphi_{i+} = (\varphi_1,\dots ,\varphi_{r_i}) \in \Gamma((q \circ f_i)^* \Sigma M)$
and each of $\varphi_1,\dots ,\varphi_{r_i} \in \Gamma(\Sigma M_i)$ is a real Killing spinor on $M_i$. 
Since $\varphi_{i+}(x) \not = 0$, at least one of $\varphi_1,\dots ,\varphi_{r_i}$ is not identically 
zero. Therefore $M_i$ admits a non-trivial real Killing spinor if $D_i \psi \not = 0$.

In a similar way \pref{11} implies that $M_i$ admits a non-trivial parallel spinor if $D_i \psi = 0$.

This completes the proof of Theorem~\ref{th1}.

\vspace{10mm}
\noindent
Bogdan Alexandrov \\
Universit\"at Greifswald \\
Institut f\"ur Mathemathik und Informatik \\
Friedrich-Ludwig-Jahn-Str. 15a \\
17487 Greifswald \\
Germany \\
{\tt e-mail: \quad boalexan@uni-greifswald.de}

\end{document}